\documentclass[12pt]{amsart}

\usepackage{latexsym}
\usepackage{amsfonts}
\usepackage{amssymb}
\usepackage{amsmath}
\usepackage{mathrsfs}
\usepackage{bbm}
\usepackage{dsfont}
\usepackage{hyperref}
\allowdisplaybreaks

\usepackage[foot]{amsaddr}

\newcommand{\be}{\begin{equation}}
	\newcommand{\ee}{\end{equation}}
\newcommand{\bea}{\begin{eqnarray}}
	\newcommand{\eea}{\end{eqnarray}}
\newcommand{\bean}{\begin{eqnarray*}}
	\newcommand{\eean}{\end{eqnarray*}}
\newcommand{\brray}{\begin{array}}
	\newcommand{\erray}{\end{array}}

\newtheorem{dfn}{Definition}[section]
\newtheorem{thm}[dfn]{Theorem}
\newtheorem{lmma}[dfn]{Lemma}
\newtheorem{ppsn}[dfn]{Proposition}
\newtheorem{crlre}[dfn]{Corollary}
\newtheorem{xmpl}[dfn]{Example}
\newtheorem{rmrk}[dfn]{Remark}

\newcommand{\bdfn}{\begin{dfn}\rm}
	\newcommand{\bthm}{\begin{thm}}
		\newcommand{\blmma}{\begin{lmma}}
			\newcommand{\bppsn}{\begin{ppsn}}
				\newcommand{\bcrlre}{\begin{crlre}}
					\newcommand{\bxmpl}{\begin{xmpl}}
						\newcommand{\brmrk}{\begin{rmrk}\rm}
							
							\newcommand{\edfn}{\end{dfn}}
						\newcommand{\ethm}{\end{thm}}
					\newcommand{\elmma}{\end{lmma}}
				\newcommand{\eppsn}{\end{ppsn}}
			\newcommand{\ecrlre}{\end{crlre}}
		\newcommand{\exmpl}{\end{xmpl}}
	\newcommand{\ermrk}{\end{rmrk}}

\newcommand{\bbc}{\mathbb{C}}
\newcommand{\bbz}{\mathbb{Z}}

\newcommand{\bbn}{\mathbb{N}}
\newcommand{\bbr}{\mathbb{R}}

\newcommand{\bbt}{\mathbb{T}}

\newcommand{\cla}{\mathcal{A}}

\newcommand{\clh}{\mathcal{H}}

\newcommand{\clk}{\mathcal{K}}
\newcommand{\cll}{\mathcal{L}}

\newcommand{\clm}{\mathcal{M}}

\title {Induced Isometric Representations}

\author{Piyasa Sarkar and S. Sundar}
\address{Institute of Mathematical Sciences (HBNI),   CIT Campus, Taramani, Chennai, India, 600113.}
\email{piyasa10@gmail.com, sundarsobers@gmail.com}

\begin{document}
	\maketitle
	
\begin{abstract}
	Let $\sigma$ be an isometric representation of $\bbn^d$ on a Hilbert space $\clh$. We induce $\sigma$ to an isometric representation $V$ of $\bbr_{+}^{d}$ on another Hilbert space $\clk$.
	We show that the map $\sigma \to V$, restricted to strongly pure isometric representations, preserves index and irreducibility. As an application, we show that, for $k \in \{0, 1,2,\cdots\} \cup \{\infty\}$, there is a continuum of prime  multiparameter CCR flows (i.e, not a tensor product
	of two non-trivial $E_0$-semigroups) with index $k$. 
	
\end{abstract}

\noindent {\bf AMS Classification No. :} {Primary 47D03 ; Secondary 46L55.}  

\section{Introduction}

Inducing representations and actions from subgroups is a time honoured method (\cite{Mackey3}, \cite{Mackey2}) to construct new representations and actions and its importance
is well established in  representation theory, ergodic theory and in many other branches of mathematics. In this paper, we consider the process of inducing isometric 
representations from subsemigroups. We only examine a  toy model where the semigroup involved is $\bbr_{+}^{d}$ and the subsemigroup involved is $\bbn^d$. 

More precisely, let $\sigma:\bbn^d \to B(\clh)$ be an isometric representation. From $\sigma$, imitating the group case, we associate an isometric representation $V$ of $\bbr_{+}^{d}$
on another Hilbert space $\clk$. We call $V$ the isometric representation induced by $\sigma$. We prove that a few properties are preserved when we pass from $\sigma$ to $V$. 
In particular, we show the following. 
\begin{enumerate}
	\item[(1)] The index of $\sigma$ coincides with the index of $V$. (For the definition of  index, see Defn. \ref{index}.)
	
	\item[(2)] Let $\{e_1,e_2,\cdots,e_d\}$ be the standard basis for $\bbn^d$. Suppose $\sigma(e_i)$ is a pure isometry for every $i=1,2,\cdots,d$. (Let us agree to call such isometric
	representations  \emph{strongly pure}.) Then, $\sigma$ is irreducible if and only if $V$ is irreducible. 
	
	\item[(3)] Let $\sigma_1,\sigma_2$ be two strongly pure irreducible isometric representations of $\bbn^d$. Denote by $V_i$ the isometric representation induced by $\sigma_i$. Then, $\sigma_1$ and $\sigma_2$ are 
	unitarily equivalent if and only if $V_1$ and $V_2$ are unitarily equivalent. 
	
\end{enumerate}
Not to get bogged down with notation, we include details only when $d=1$ or $d=2$. 

The takeaway from $(2)$ and $(3)$ is that 'enumerating' irreducible isometric representations of $\bbr_{+}^{d}$ is at least as hard as enumerating irreducible isometric representations of $\bbn^d$.
For $d=2$, it is known from \cite{Berger} that irreducible isometric representations of $\bbn^2$, except the one-dimensional ones, are in one-one correspondence with the irreducible unitary representations of $\bbz_2\ast \bbz$ whose associated group $C^{*}$-algebra is not type I, and consequently its representation theory is quite complicated. Thus, unlike the $1$-parameter case, the classification problem of  isometric representations of $\bbr_{+}^{d}$, when $d \geq 2$, is quite hard. 

Why consider induced isometric representations? The motivation for us come from a problem in the multiparameter theory of $E_0$-semigroups which we next explain. 
Let us recall a few definitions regarding $E_0$-semigroups. Let $P$ be a closed, convex cone in $\bbr^d$ that spans $\bbr^d$ and contains no line. Let $\clh$ be a separable
Hilbert space. By an $E_0$-semigroup, over $P$, on $B(\clh)$, we mean a semigroup $\alpha:=\{\alpha_x\}_{x \in P}$ of unital normal $*$-endomorphisms of $B(\clh)$ such that 
the map 
\[
P \ni x \to \langle \alpha_x(A)\xi|\eta \rangle \in \bbc\]
is continuous for every $A \in B(H)$, and $\xi,\eta \in \clh$. The equivalence relation on `the set of $E_0$-semigroups' that we consider is that of cocycle conjugacy. The first numerical invariant  for $E_0$-semigroups, in the $1$-parameter case, is  due to Arveson (\cite{Arv_Fock}, \cite{Arveson}), and is called the index. Arveson proved that index is a complete invariant for $1$-parameter CCR flows.  Moreover, index is only relevant for spatial $E_0$-semigroups.   

An $E_0$-semigroup $\alpha:=\{\alpha_x\}_{x \in P}$ on $B(\clh)$ is said to be \emph{spatial} if it has a unit, by which we mean a strongly continuous semigroup of bounded operators $u:=\{u_x\}_{x \in P}$ on $\clh$ such that 
\begin{enumerate}
	\item[(1)] for $x \in P$, $u_x \neq 0$, and
	\item[(2)] for $A \in B(\clh)$, $x \in P$, $\alpha_x(A)u_x=u_xA$.
\end{enumerate} 
Following  Arveson, in \cite{sarkar_sundar_2022} the authors, in the multiparameter case,  defined a numerical invariant for a spatial $E_0$-semigroup $\alpha$ called the index of $\alpha$ which we denote by $Ind(\alpha)$. Roughly, $Ind(\alpha)$ measures `the number of units' of $\alpha$. 

Arveson, in the $1$-parameter case, proved the remarkable fact that index is a `homomorphism', i.e. 
\begin{equation}
	\label{additivity}
	Ind(\alpha \otimes \beta)=Ind(\alpha)+Ind(\beta).\end{equation}
We mention here that Arveson's proof, without any modification, works in the multiparameter case as well. In view of the above equation, it is quite natural to ask the following question. 

\textbf{Question:} Does there exist a prime $E_0$-semigroup that has index at least two?

Recall that an $E_0$-semigroup $\alpha$ is said to be \emph{prime} if whenever $\alpha$ is cocycle conjugate to $\beta \otimes \gamma$, where $\beta$ and $\gamma$ are $E_0$-semigroups,  then either $\beta$ or $\gamma$ is an automorphism group. In short, an $E_0$-semigroup is prime if it cannot be written as a tensor product of two non-trivial $E_0$-semigroups. 
Affirmative answer to the above question, when $P=[0,\infty)$, is quite complicated and was given by Liebscher (\cite{Liebscher}) who constructed such examples by probabilistic 
means. However, in the multiparameter case we show that such examples exist even within the class of CCR flows (which are probably the first examples studied in the theory of $E_0$-semigroups);  a  total contrast to the one parameter case. 

It follows from the results of \cite{sarkar_sundar_2022} and \cite{SUNDAR2021108802}, which are collected in Section 2, that the problem of constructing prime CCR flows over $P$ with a given index $k$  is equivalent to 
the problem of constructing irreducible isometric representations of $P$ with index $k$. 
Here is where induced isometric representations come into picture. For the discrete semigroup  $\bbn^2$, such examples are available in the literature (\cite{ALBEVERIO201935}), albeit in a slightly disguised form.  We also construct alternate examples. We induce such discrete semigroups of isometries to construct the desired isometric representations, and we prove the following. 

\begin{thm}\label{main_intro}
	Let $P$ be a closed convex cone in $\bbr^d$ which is spanning and pointed. Suppose that $d \geq 2$. Then, for each $k \in \{0,1,2,\cdots,\} \cup \{\infty\}$, there is a continuum of irreducible isometric representations of $P$ that has index $k$. 
\end{thm}

The following  theorem is now immediate. 

\begin{thm}
\label{main prime}
	Let $P$ be a closed convex cone in $\bbr^d$ which is pointed and spanning. Suppose that $d \geq 2$. Then, for each $k \in \{0, 1,2,\cdots,\} \cup \{\infty\}$, there is a continuum of prime CCR flows with index $k$. 
\end{thm}

We end this introduction by mentioning that for $k \in \{0,1\}$, the above theorem is known. For $k=0$, the CCR flows considered in  \cite{Anbu_Sundar}  provide such examples. For $k=1$, the authors in \cite{sarkar_sundar_2022} constructed such examples. 
We must mention here that the examples constructed in \cite{sarkar_sundar_2022} are the first `genuine' examples of  CCR flows/$E_0$-semigroups, in the multiparameter case, that are type one (which roughly means that there is abundance of units). In \cite{sarkar_sundar_2022}, the focus was on to construct type one examples with index one. On taking tensor product of such  examples, we can easily construct type one examples with index greater than one. But this is clearly tautological and this motivated us to seek examples of prime CCR flows with index greater than one. 

\textit{Notation:-} For us, $\bbn$ stands for the set of  natural numbers together with $0$. We denote $[0,\infty)$ by $\bbr_{+}$. Our convention is that inner products are linear in the first variable.

\section{Preliminaries}
	
First, we recall a few definitions that we need. 		Let $G$ be a locally compact, abelian, second countable, Hausdorfff topological group, and let $P \subseteq G$ be a closed semigroup containing $0$ such that $P-P =G$ and $\overline{Int(P)} = P$. Let $\clh$ be a separable Hilbert space. 
	Let $V=\{V_x\}_{x\in P}$ be a strongly continuous semigroup of isometries on $\clh$. Such a family is also called an \textsl{isometric representation} of $P$ on $\clh$. We call $V$ a pure isometric representation of $P$ if $\displaystyle \bigcap_{x \in P} Ran (V_x) = \{0\}$. The representation $V$ is said to be \textsl{irreducible} if the only closed subspaces of $\clh$ invariant under $\{V_x, V_x^* | x \in P \}$ are $\{0\}$ and $\clh$.
	
	Let $V= \{V_x\}_{x\in P}$ be an isometric representation of $P$ on a Hilbert space $\clh$. A   map $\xi :P \to \clh$ is called an \textsl{additive cocycle of V} if
	\begin{itemize}
		\item [(a)] for all $x$ in $P$, $\xi_x \in ker(V_x^*)$, and
		\item[(b)] for all $x,y$ in $P$, $\xi_{x+y} = \xi_x + V_x \xi_y$.
	\end{itemize}
	The vector space of all additive cocycles of $V$ is denoted by $\mathcal{A}(V)$. 
	
	\begin{rmrk}
	Let $V:P \to B(\clh)$ be an isometric representation, and let $\xi:P \to \clh$ be an additive cocycle. Then, $\xi$ is norm continuous. To see this, let $(s_n)_{n \geq 1}$ be a cofinal sequence in $Int(P)$. Set $E_n:=\{x \in P: s_n-x \in Int(P)\}$. Note that for $x \in E_n$, 
	\[
	\xi_x=(1-V_xV_{x}^{*})\xi_{s_n}.\]
	It follows from the above equality and the fact that $(E_n)_{n \geq 1}$ is an open cover of $P$ that $\xi$ is norm continuous. 
	
	\end{rmrk}
	
	\begin{dfn} \label{index}
	 For an isometric representation $V$ of $P$, we define the \textsl{index of V}, denoted $Index(V)$, as the dimension of $\mathcal{A}(V)$.
	 \end{dfn} 
 Note that if $V=V_1 \oplus V_2$, then. 
	\[ Index(V) = Index(V_1) +Index(V_2).\]	
	\textit{Notation:-} We define $\clm(V) := C^* (\{ V_x, V_x^*| x \in P\})$.	
	
	\begin{ppsn} \label{smaller cone}
		Let $G$ be a locally compact abelian group, and let $P$ be a closed semigroup of $G$ containing $0$ such that $P-P=G$. Let $V:P \to B(\clh)$ be a pure isometric representation on a separable Hilbert space $\clh$. Let $Q$ be another closed semigroup of $G$ containing $0$ such that $Q-Q=G$ and $Q \subset P$. Denote the restriction of $V$ to $Q$ by $W$. Then,
		\begin{enumerate}
			\item [(1)] $W$ is pure,
			\item [(2)] $dim(\cla(W)) = dim(\cla(V))$, and
			\item [(3)] $\clm(W)' = \clm(V)'$.
		\end{enumerate}
	Suppose $V^{(1)}$ and $V^{(2)}$ are two pure isometric representations of $P$ acting on Hilbert spaces $\clh_1$ and $\clh_2$ respectively. Denote the restrictions of $V^{(i)}$ to $Q$ by $W^{(i)}$. Then, $V^{(1)}$ and $V^{(2)}$ are unitarily equivalent if and only if $W^{(1)}$ and $W^{(2)}$ are unitarily equivalent. 
	\end{ppsn}
	
	\textit{Proof.} To see that $W$ is a pure isometric representaion, assume $\xi \in \bigcap_{x \in Q} Ran(W_x)$. Let $a \in P$. Since $Q$ is spanning in $G$, there exist $x,y \in Q$ such that $a=x-y$, i.e $a+y=x$. Now, there exists $\eta$ such that $\xi=W_x \eta = V_x \eta = V_a(V_y \eta)$. This implies $\xi \in Ran(V_a)$, for all $a \in P$. However $V$ is pure, and thus $\xi=0$. Hence, $W$ is pure.
	
	Let $\xi = \{ \xi_x\}_{x \in P}$ be an additive cocycle of $V$. Define, for $ y \in Q$,
	\[\eta^{\xi}_y := \xi_y.\]
	It is straightforward to see that $\eta^{\xi} = \{ \eta^{\xi}_y\}_{y \in Q} \in \cla(W)$. We claim that the map 
	\[ \cla(V) \ni \xi \mapsto \eta^{\xi} \in \cla(W) \]
	is an isomorphism. To see that it is injective, suppose $\xi \in \mathcal{A}(V)$ is such that $\xi_y=0$ for all $y \in Q$. Then, for $a \in P$, write $a=x-y$ with $x,y \in Q$, and calculate to observe that 
	\[
	0=\xi_x=\xi_{a+y}=\xi_a+V_a\xi_y=\xi_a.\]
	Thus, $\xi=0$, and hence  the map is injective. 
	
	For proving it is a surjection, let $\eta$ be an additive cocycle of $W$. Note that for any $c,d,\alpha \in Q$, and $b \in P$ such that $b=c-d$,
	\begin{align*}
		\eta_{c + \alpha} - V_b \eta_{d+ \alpha} &= \eta_c + W_c \eta_{\alpha} - V_b \eta_d - V_b W_d \eta_{\alpha} \\
		&= \eta_c + V_c \eta_{\alpha} - V_b \eta_d - V_c \eta_{\alpha} \\
		&= \eta_c - V_b \eta_d.
	\end{align*}
	Thus, for $c,d,\alpha \in Q$, and $b \in P$, if $b=c-d$, then
	\begin{equation}
	\label{well defined nature}
			\eta_{c + \alpha} - V_b \eta_{d+ \alpha}=\eta_c - V_b \eta_d.\end{equation}
	Let $a \in P$. Then, there exist $x, y \in Q$ such that $a=x-y$. Define 
	\[ \xi_a :=  \eta_x - V_a \eta_y.\]
	Say there also exist $u,v \in Q$ such that $a = u-v$. This implies $x+v = y+u$, and applying Eq. \ref{well defined nature} twice, we get 
	\[\eta_x - V_a \eta_y = \eta_{x+v} - V_a \eta_{y+v} = \eta_{u+y} - V_a \eta_{v+y} = \eta_u - V_a \eta_v.\] 
	Thus, $\xi_a$ is well-defined. Also, 
	\begin{align*}
		V_a^* \xi_a &= V_a^*\eta_x -  \eta_y \\
		&= V_x^* V_y \eta_x -\eta_y = V_x^* (\eta_{x+y} - \eta_y) - \eta_y \\
		&=V_{x}^{*}(\eta_x+V_x\eta_y-\eta_y)-\eta_y\\
		&= \eta_y - V_x^* \eta_y - \eta_y \\
		&= -V_x^* \eta_y \\
		&= -V_a^* V_y^* \eta_y = 0.
	\end{align*}
	This shows that $\xi_a \in ker(V_a^*)$. To prove the cocycle nature, let $a, b \in P$. There exist $x,y,z$ such that $a=x-y$ and $b=y-z$. Then, 
	\begin{align*}
		\xi_{a+b} &= \eta_x - V_{a+b} \eta_z \\
		&=\eta_x-V_a\eta_y+V_a\eta_y-V_{a+b}\eta_z\\
		&=\xi_a+V_a(\eta_y-V_b\eta_z)\\
		&=\xi_a+V_a\xi_b. 
	\end{align*}
	Therefore, $\xi = \{ \xi_x\}_{x \in P} \in \cla(V)$. Also, for $y \in Q$, note that $\xi_y = \eta_y$. Thus, $\eta^{\xi} = \eta$, and the map is hence a bijection.
	
	Note that for $a \in P$ if $a=x-y$ with $x,y \in Q$, $V_{y}V_a=V_x$ which in turn  implies implies that  $V_a=V_{y}^{*}V_x=W_{y}^{*}W_x$. Thus, the  $C^*$-algebras generated by $\{V_a| a \in P\}$ and $\{W_x| x \in Q\}$ are the same and that implies $\clm(W)' = \clm(V)'$.
	
	Let $V^{(1)}$ and $V^{(2)}$ be two pure isometric representations of $P$ acting on Hilbert spaces $\clh_1$ and $\clh_2$ respectively. Denote the restrictions of $V^{(i)}$ to $Q$ by $W^{(i)}$. If $V^{(1)}$ is unitarily equivalent to $V^{(2)}$, it is clear that $W^{(1)}$ is uniatrily equivalent to $W^{(2)}$. Conversely, if $W^{(1)}$ is unitarily equivalent to $W^{(2)}$, there exists a unitary $U: \clh_1 \to \clh_2$ such that $W^{(1)}_x U = U W^{(2)}_x$, for all $x \in Q$. Let $a \in P$; there exist $x, y \in Q$ such that $a=x-y$. Then, 
	\begin{align*}
		V^{(1)}_a U &= W^{(1)*}_y W^{(1)}_x U \\
		&= W^{(1)*}_y U W^{(2)}_x \\
		&=   U W^{(2)*}_y W^{(2)}_x \\
		&= U V^{(2)}_a.
	\end{align*}
	This proves that $V^{(1)}$ is unitarily equivalent to $V^{(2)}$ iff $W^{(1)}$ is unitarily equivalent to $W^{(2)}$.  \hfill $\Box$

	Now, let $P$ be a closed, convex one in $\bbr^d$ that spans $\bbr^d$ and is pointed. 
	The simplest class of examples of $E_0$-semigroups over $P$ are the CCR flows which  arise from isometric representations of $P$. Let us recall the definition of the CCR flow associated to an isometric representation. 
	
	Let $V= \{V_x \}_{x\in P}$ be an isometric representation of $P$ on a Hilbert space $\clh$. Let $\Gamma(\clh)$ denote the symmetric Fock space of $\clh$. There exists a unique $E_0$-semigroup, denoted $\alpha^V$, on $B(\Gamma(\clh))$, such that for all $x \in P$ and $\xi \in \clh$,
	\[\alpha^V_x(W(\xi)) = W(V_x \xi) \]
	where $\{W(\xi)| \xi \in \clh \}$ is the collection of Weyl operators on $\Gamma(\clh)$. We call $\alpha^V$ the \textsl{CCR flow associated to V}. 

	\begin{rmrk}
	\label{on ccr}
	We collect a few facts concerning CCR flows in this remark. 
	\begin{enumerate}
	\item[(1)] For two pure isometric representations $V_1$ and $V_2$, the CCR flows $\alpha^{V_1}$ and $\alpha^{V_2}$ are cocycle conjugate if and only if $V_1$ and $V_2$ are unitarily equivalent. For a proof, we refer the reader to Thm. 5.2 of \cite{SUNDAR2021108802}
	\item[(2)] Let $V=\{V_x\}_{x \in P}$ be a pure isometric representation of $P$ on a Hilbert space $\clh$. 
	\begin{enumerate}
		\item[(a)] It follows from Thm. 7.2 of \cite{SUNDAR2021108802} that the CCR flow $\alpha^V$ is prime if and only if 
		 the representation $V$ is irreducible.
		\item[(b)] By Prop. 2.7 of  \cite{sarkar_sundar_2022}, it follows that $Ind(\alpha^V)=dim (\cla (V))$.
 			\end{enumerate}
	\end{enumerate}

	\end{rmrk}

		\section{Induced isometric Representations}
	Let $\clh$ be a separable Hilbert space with an orthonormal basis $\{ e_n\}_{n \in \mathbb{N}}$. Let $d \geq 1$ be an integer. Let $\sigma: \mathbb{N}^d \to B(\clh)$ be an isometric representation. For any $(x_1,...x_d)$ in $[0, \infty)^d$, we denote it by $\tilde{x}$, and elements in $\mathbb{N}^d$ by $\tilde{n}$. Let
	\begin{align*}
		\clk :=& \{ \xi: [0,\infty)^d \to \clh |\xi \textrm{ is measurable, square-integrable over compact sets and } \\ & \xi( \tilde{x} + \tilde{n}) = \sigma(\tilde{n}) \xi(\tilde{x}), \forall \tilde{x} \in [0,\infty)^d, \tilde{n} \in \mathbb{N}^d \}
	\end{align*} 
	Define an inner product $\langle ~ | ~ \rangle $ on $\clk$ by
	\[ \langle \xi| \eta \rangle := \int_{0}^{1} \int_{0}^{1}... \int_{0}^{1} \langle \xi(\tilde{x}) | \eta(\tilde{x}) \rangle d(\tilde{x}) \]
	for all $\xi, \eta \in \clk$. It goes without saying that we identify two elements of $\mathcal{K}$ if they agree almost everywhere. 
	Then, $\clk$ is a Hilbert space under this inner product. For each $ \tilde{t} \in [0,\infty)^d$, define $ V_{\tilde{t}}: \clk \to \clk $ by
	\[ V_{\tilde{t}} \xi (\tilde{x}) := \xi ( \tilde{x} + \tilde{t}).\]
	It is routine to check that $V:=\{V_{\tilde{t}}\}_{\tilde{t} \in \bbr_+^d}$ is a strongly continuous semigroup of isometries. The representation $V: [0, \infty)^d \to B(\clk)$ thus obtained is an isometric representation, and we shall call it the \emph{isometric representation induced by $\sigma$}.

	The space $\clk$ can be identified with $L^2([0,1)^d) \bigotimes \clh$ via the map
	\[ \clk \ni \xi \mapsto \sum_{n \in \mathbb{N}} (\xi_n \otimes e_n) \in L^2([0,1)^d) \bigotimes \clh, \]
	where $\xi_n: [0,1)^d \to \mathbb{C}$ is defined by $\xi_n(\tilde{x}) = \langle \xi(\tilde{x})| e_n \rangle$, for all $\tilde{x} \in [0,1)^d, n \in \mathbb{N}$. We always use this identification.  Under this identification, we obtain \[V_{\tilde{n}} = 1_{L^2([0,1)^d)} \otimes \sigma(\tilde{n}), \textrm{for all }\tilde{n} \in \mathbb{N}^d.\]
	
	\begin{rmrk}
		\begin{enumerate}
			\item 	Note that if $\sigma$ is a pure isometric representation on $\clh$, then $V$ too is a pure isometric representation on $\clk$. To see this, consider a pure isometric representation $\sigma: \mathbb{N}^d \to B(\clh)$ and the corresponding induced representation $V: [0,\infty)^d \to B(\clk)$. This implies, for all $\tilde{n} \in \mathbb{N}^d$,  $V_{\tilde{n}} = 1_{L^2([0, \infty)^d)} \otimes \sigma(\tilde{n})$. Since $\sigma$ is pure, $\displaystyle \bigcap_{\tilde{n} \in \mathbb{N}^d} Ran(\sigma(\tilde{n})) = \{0\}$. Thus, \[\displaystyle  \bigcap_{\tilde{t} \in [0, \infty)^d} Ran(V_{\tilde{t}}) \subseteq \bigcap_{\tilde{n} \in \mathbb{N}^d} Ran(V_{\tilde{n}}) = \{0\}.\] Hence, $V$ is a pure isometric representation.
			
			\item Suppose $\sigma^{(1)}$ and $\sigma^{(2)}$ are two pure irreducible isometric representations of $\bbn^d$ acting on $\clh_1$ and $\clh_2$ respectively, and let $V^{(1)}$ and $V^{(2)}$ be the corresponding induced representations. Then, $\sigma^{(1)}$ and $\sigma^{(2)}$ are unitarily equivalent if and only if $V^{(1)}$ and $V^{(2)}$ are unitarily equivalent. This follows from Schur's lemma and the fact that for $i=1,2$, and $\tilde{n} \in \bbn^d$, $V^{(i)}_{\tilde{n}}=1 \otimes \sigma^{(i)}_{\tilde{n}}$.
		
		\end{enumerate}

	\end{rmrk}

	We now look at induced representations when $d=1$ and $d=2$.
	
	Let us consider the case when $d=1$. Let $\clh$ be a separable Hilbert space with an orthonormal basis $\{e_n\}_{n \in \mathbb{N}}$. Let $\sigma: \mathbb{N} \to B(\clh)$ be a pure isometric representation on $\clh$. Let $V$ be the isometric representation induced by $\sigma$ on $\clk$. Recall that
	\begin{align*}
		\clk = \{ & \xi: [0,\infty) \to \clh | \xi \textrm{ is measurable, square-integrable over compact sets, }\\ &\xi(t+n)=\sigma(n)\xi(t), \forall t \geq 0, n \in \mathbb{N} \},
	\end{align*}
	where the inner product is defined by 
	\[ \langle \xi| \eta \rangle = \int_{0}^{1} \langle \xi(t)|\eta(t) \rangle dt. \]
	We usually  identify $\clk$ with $ L^2([0,1))\bigotimes \clh$ and  `move' between the two spaces freely whenever convenient. Recall that, for each $t\geq0$, $ V_t: \clk \to \clk$ is given by,
	\[ V_t \xi (x) := \xi(x +t).\]
	The representation $V= \{ V_t\}_{t \geq 0} $ is pure, since $\sigma$ is pure.
		We first calculate $V_t^*$. Let $t \geq 0$, there exists $n \in \mathbb{N}$ such that $n \leq t < n+1$. Let $\xi, \eta \in \clk$. Then,
	\begin{align*}
		~~~\langle V_t^* \xi | \eta \rangle &= \langle \xi | V_t \eta \rangle \\
		&= \int_{0}^{1} \langle \xi(x) | \eta(x+t) \rangle dx \\
		&= \int_{0}^{1} \langle \sigma(n)^* \xi(x) | \eta(x+t-n) \rangle dx = \int_{t-n}^{t-n+1} \langle \sigma(n)^* \xi(x-t+n) | \eta(x) \rangle dx \\
		&= \int_{t-n}^{1} \langle \sigma(n)^* \xi(x-t+n) | \eta (x) \rangle dx + \int_{1}^{1+t-n} \langle \sigma(n)^* \xi(x-t+n)| \eta(x) \rangle dx \\
		&= \int_{t-n}^{1} \langle \sigma(n)^* \xi(x-t+n) | \eta(x) \rangle dx + \int_{1}^{1+t-n} \langle \sigma(n+1)^* \xi (x-t+n) | \eta (x-1) \rangle dx \\
		&= \int_{t-n}^{1} \langle \sigma(n)^* \xi(x-t+n) | \eta (x) \rangle dx + \int_{0}^{t-n} \langle \sigma(n+1)^* \xi (x-t+n+1)| \eta (x) \rangle dx
	\end{align*}
	Thus, we have  
	\begin{align}\label{V_t*}
		V_t^* \xi(x) = \begin{cases}
			\sigma(n+1)^* \xi (x-t+n+1), \textrm{ for } x \in [0,t-n), \\
			\sigma(n)^* \xi (x-t+n), \textrm{ for } x \in [t-n,1).
		\end{cases}
	\end{align}
	Let $ 0 \leq t <1$. Note that Eqn. \ref{V_t*} implies that
	\begin{align} \label{ker(V_t^*)}
		ker(V_t^*) = \{ \xi \in L^2((0,1], \clh) ~|~ \xi|_{[0,1-t)} = 0 \textrm{ and } \xi(x) \in ker(\sigma(1)^*), ~ x \in [1-t,1)\}
	\end{align}
	Also, keeping in mind the correspondence $\clk \cong L^2((0,1], \clh)$, for $n \leq t < n+1$ and $\xi \in \clk$,
	\begin{align}\label{V_t}
		(V_t \xi)(x) = \begin{cases}
			\sigma(n)\xi(x+t-n), ~ x \in [0,1-t+n), \\
			\sigma(n+1) \xi(x+t-n-1), ~ x \in [1-t+n,1)
		\end{cases}
	\end{align}
	
	\begin{lmma} \label{iff}
		Let $\eta: \mathbb{N} \to \clh$ be a map. Define, for $n \leq t < n+1$, $\xi^{\eta}_t : [0,1) \to \clh$,
		\[  \xi^{\eta}_t(x) := \begin{cases} \eta_n, ~ x \in [0, 1-t+n), \\
			\eta_{n+1}, ~ x \in [1-t+n, 1).
		\end{cases} \]
		Then, $\xi^{\eta} = \{\xi^{\eta}_t\}_{t \geq 0}$ is an additive cocycle of $V$ iff $\eta = \{\eta_k\}_{k \in \mathbb{N}}$ is an additive cocycle of $\sigma$.
	\end{lmma} 
	
	\textit{Proof.} 
	 Assume that $\xi^{\eta}$ is an additive cocycle of $V$. Since $\xi^{\eta}_k \in ker(V_k^*)$, 
	\[ \sigma(k)^* \eta_k= \sigma(k)^* \xi^{\eta}_k(x) =V_k^* \xi^{\eta}_k (x)= 0.\]
	Thus, $\eta_k \in ker(\sigma(k)^*)$, for all $k \in \mathbb{N}$. Also, 
	\[\eta_{k+m}= \xi^{\eta}_{k+m}(x )=  \xi^{\eta}_k(x) + V_k\xi^{\eta}_m(x) = \eta_k + \sigma(k) \eta_m.\]
	Thus, $\eta = \{ \eta_k\}_{k \in \mathbb{N}}$ is an additive cocycle of $\sigma$.
	
	Suppose $\eta $ is an additive cocycle of $\sigma$. Then, for $n \leq t < n+1$,
	\begin{align*}
		V_t^* \xi^{\eta}_t (x) &= \begin{cases} \sigma(n+1)^* \xi^{\eta}_t(x-t+n+1), \textrm{ whenever } x \in [0, t-n), \\ \sigma(n)^* \xi^{\eta}_t (x-t+n), \textrm{ whenever } x \in [t-n,1). \end{cases} \\
		&= \begin{cases} \sigma(n+1)^* \eta_{n+1}, \textrm{ whenever } x \in [0, t-n), \\ \sigma(n)^* \eta_n , \textrm{ whenever } x \in [t-n,1). \end{cases} \\
		&=0.
	\end{align*} 
	This shows $\xi^{\eta}_t \in ker(V_t^*)$, for all $t \geq 0$.

	 To prove the cocycle nature, let $s,t \geq 0 $ be given. Choose $m,n \in \mathbb{N}$ such that $ m \leq s < m+1$ and $n \leq t < n+1$.Then,
        \begin{enumerate}
                \item [Case (i)] $m+n \leq s+t < m+n+1$, \\
                Under this condition,
                \[ \xi^{\eta}_{s+t}(x) = \begin{cases}
                        \eta_{m+n}, ~ x \in [0,1-s-t+m+n),\\
                        \eta_{m+n+1}, ~ x \in [1-s-t+m+n,1).
                \end{cases}\]
                Now,
                \begin{align*}
                        \xi^{\eta}_s(x) + V_s \xi^{\eta}_t(x) &= \begin{cases}
                                \xi^{\eta}_s(x) + \sigma(m) \xi^{\eta}_t(x+s-m), ~ x \in [0,1-s+m), \\
                                \xi^{\eta}_s(x) + \sigma(m+1) \xi^{\eta}_t(x+s-m-1), ~ x \in [1-s+m,1). \\
                        \end{cases} \\
                        &= \begin{cases}
                                \eta_m + \sigma(m) \eta_n, ~ x \in [0,1-s+m) \cap [0,1-s-t+m+n), \\
                                \eta_m + \sigma(m) \eta_{n+1}, ~ x \in [0,1-s+m) \cap [1-s-t+m+n,1), \\
                                \eta_{m+1} + \sigma(m+1) \eta_n, ~ x \in [1-s+m,1) \cap [0,1). \\
                        \end{cases} \\
                        &= \begin{cases}
                                \eta_{m+n}, ~ x \in [0,1-s-t+m+n),\\
                                \eta_{m+n+1}, ~ x \in [1-s-t+m+n,1).\\
                        \end{cases} \\
                        &= \xi^{\eta}_{s+t}(x).
                \end{align*}

                \item [Case (ii)] $m+n+1 \leq s+t < m+n+2$. The verification in this case is similar. 
        \end{enumerate}
        Thus, $\xi^{\eta}=\{ \xi^{\eta}_t\}_{t \geq 0}$ is an additive cocycle of $V$ whenever $\eta$ is an additive cocycle of $\sigma$. 
	
	\hfill $\Box$

	Let $\cla(\sigma)$ and $\cla(V)$ denote the space of additive cocycles of $\sigma$ and $V$ respectively. 
	
	\begin{ppsn} \label{1-d bijection}
		The map 
		\[ \cla(\sigma) \ni \eta \mapsto \xi^{\eta} \in \cla(V)\]
		is an isomorphism. Hence, $Ind(\sigma)=Ind(V)$.
	\end{ppsn}
	
	\textit{Proof.} The map is clearly injective. To see that it is also onto, let $\xi= \{ \xi_t\}_{t \geq 0}$ be a non-zero additive cocycle of $V$. Recall that $V$ acts on $\clk$ where $\clk$ is given by 
	\begin{align*}
		\clk = \{ & \xi: [0,\infty) \to \clh | \xi \textrm{ is measurable, square-integrable over compact sets, }\\ &\xi(t+n)=\sigma(n)\xi(t), \forall t \geq 0, n \in \mathbb{N} \},
	\end{align*}
	The fact that $\xi$ is an additive cocylce implies that  for all $s,t \geq 0$, $\xi_{s+t}(x)= \xi_s(x)+\xi_t(x+s)$ a.e. By Theorem 5.3.2 of \cite{Arveson}, there exists  $f \in L^2_{loc}([0,\infty), \clh)$, such that for every $t$,
	\begin{equation} \label{Arveson trick} \xi_t(x)= f(x+t) - f(x),\end{equation}
	for almost all $x$. However, $\xi_t \in ker(V_t^*)$, and hence, by Eqn. \ref{ker(V_t^*)}, for $0 \leq t <1$, $\xi_t(x) = 0$, for almost all $x \in [0,1-t)$. This implies, for $0 \leq t <1$, $f(x+t)=f(x)$, for almost all $x \in [0,1-t)$. Define $f_n: [0, \infty) \to \mathbb{C}, f_n(x):= \langle f(x)|e_n \rangle$, for all $n \in \mathbb{N}$. Then, for every $t \in [0,1)$, $f_n(x+t)=f_n(x)$ for almost all $x \in [0,1-t)$. Thus, the distributional derivative of $f_n$ is zero on $(0,1)$. Hence, the function $f_n$ is constant on the interval $[0, 1)$, and thereby, so is the function $f$
		
	  Therefore, there exists a vector $\gamma \in \clh$, such that \begin{equation}
	  \label{value of f}
	  f(x)= \gamma \end{equation} for almost all $x \in [0,1)$. 
	  
	Let $\phi:[0,\infty) \to \clh$ be defined by \[\phi(x)=f(x+1)-\sigma(1)f(x).\]  For every $t>0$, since $\xi_{t}(x+1)=\sigma(1)\xi_t(x)$ for almost all $x$, it follows that given $t>0$, \[f(x+t+1)-f(x+1)=\sigma(1)(f(x+t)-f(x))\]for almost all $x \in [0,\infty)$. In other words, given $t>0$, $\phi(x+t)=\phi(x)$ for almost all $x \in [0,\infty)$. This forces that $\phi$ is constant. Let  $\widetilde{\gamma} \in \clh$ be such that $\phi(x)=\widetilde{\gamma}$ for almost all $x \in [0,\infty)$. It follows from Eq. \ref{value of f}, that for almost all $x \in [0,1)$, 
	\begin{equation}
	\label{value of f beyond one}
	f(x+1)=\sigma(1)f(x)+\widetilde{\gamma}=\sigma(1)\gamma+\widetilde{\gamma}
	\end{equation}
	for almost all $x \in [0,1)$. Set $\eta_1:=\sigma(1)\gamma+\widetilde{\gamma}-\gamma$. 
	  
	  Combining Eq. \ref{Arveson trick}, Eq. \ref{value of f} and Eq. \ref{value of f beyond one},  we have, for $0 \leq t <1$, 
	\[ \xi_t(x) = \begin{cases} 0 ,~ x \in [0,1-t), \\
		\eta_1, ~ x \in [1-t,1). \end{cases} \]
	
	By Eqn. \ref{ker(V_t^*)}, we have $\eta_1 \in ker(\sigma(1)^*)$. 	Set $\eta_{k+1} = \eta_k + \sigma(k) \eta_1$, for all $k \geq 1, k \in \mathbb{N}$, and $\eta_0=0$. This makes $\eta= \{ \eta_k\}_{k \in \mathbb{N}}$ an additive cocycle of $\sigma$. Therefore, $\xi^{\eta}$ is an additive cocycle of $V$. 
	Also, \[\xi_t=\xi^{\eta}_t\] for $0 \leq t <1$. Since $\xi$ and $\xi^\eta$ are additive cocycles, $\xi=\xi^\eta$. 
	  This proves that the map is indeed a bijection. \hfill $\Box$

	\begin{rmrk}
		Let $W= \{W_t\}_{t \geq 0}$ be a pure isometric representation on a separable Hilbert space $\cll$. Then, $\{W_t\}_{t \geq 0}$ is unitarily equivalent to $\{S_t \otimes 1\}_{t \geq0 }$ on $L^2([0,\infty)) \bigotimes \clh_0$, for some separable Hilbert space $\clh_0$. Here, $\{S_t\}_{t \ge 0}$ is the shift semigroup on $L^2([0,\infty))$. Moreover,  $dim(\cla(W))= dim(\clh_0)$. The last equality can be proved directly, or by appealing to the index computation, due to Arveson, of $1$-parameter CCR flows and the fact that for CCR flows, $Ind(\alpha^W)=dim(\cla(W))$. Thus, $W$ is irreducible if and only if  $dim(\cla(W))=1$.	
	\end{rmrk}
	
	Recall that
	\begin{align*}
		\clm(\sigma) = C^*(\{ \sigma(1)\}), \\
		\clm(V) = C^*( \{ V_t | t \geq 0 \})
	\end{align*}
	and denote their commutants by $\clm(\sigma)'$ and $\clm(V)'$.

	\begin{ppsn} \label{1-d commutant}
		With the foregoing notation, $\clm(V)' = 1_{L^2([0,1))} \bigotimes \clm(\sigma)'. $	
	\end{ppsn}
	
	\textit{Proof.} Define $\tilde{\sigma} : \mathbb{N} \to B(\ell^2(\mathbb{N}))$, $\tilde{\sigma}(1) = S$, where $S: \ell^2(\mathbb{N}) \to \ell^2(\mathbb{N})$ is the unilateral shift operator. Then, $\tilde{\sigma}$ is a pure isometric representation on $\ell^2(\mathbb{N})$. Let $\tilde{V}: [0, \infty) \to B(L^2([0,1)) \bigotimes \ell^2(\mathbb{N}))$ be the pure isometric representation induced by $\tilde{\sigma}$. Clearly,  $dim(\cla(\tilde{\sigma}))=1$. By Prop. \ref{1-d bijection}, we have $dim(\cla(\tilde{V}))=1$, and consequently, $\tilde{V}$ is irreducible, and
	\begin{equation}
		\label{commutant1}
		\clm(\tilde{V})' = \mathbb{C} 1_{L^2([0,1))} \otimes 1_{\ell^2(\mathbb{N})}.
	\end{equation}

	Let $\sigma:\bbn \to B(\clh)$ be a pure isometric representation. Since $\sigma: \mathbb{N} \to B(\clh)$ is  pure, by Wold Decomposition, there exists a Hilbert space $\clh_0$ and a unitary $U: \ell^2(\mathbb{N}) \bigotimes \clh_0 \to \clh$ such that 
	\[\sigma(1)= U(S \otimes 1) U^*.\]
	So, we can assume that up to a unitary equivalence, $\clh = \ell^2(\mathbb{N}) \otimes \clh_0$, for some Hilbert space $\clh_0$, and $\sigma(1)=S \otimes 1=\tilde{\sigma}(1) \otimes 1$,  where $1 \in B(\clh_0)$ is the identity on $\clh_0$. Let $V$ be the induced representation due to $\sigma$.   We thus have 
	\begin{align*}
		\clm(\sigma)' &= \{ 1_{\ell^2(\mathbb{N})} \otimes T| T \in B(\clh_0)\}, \\
		V_t &= \tilde{V_t} \otimes 1, ~\textrm{ and} \\
		\clm(V)' &= \{ 1_{L^2([0,1))} \otimes 1_{\ell^2(\mathbb{N})} \otimes T | T \in B(\clh_0)\} ~~(\textrm{by Eq. \ref{commutant1}}) \\
		&= 1_{L^2([0,1))} \otimes \clm(\sigma)'. 
	\end{align*}
	The proof is complete. 	 \hfill $\Box$

	Let us now consider the case when $d=2$. 	Let $\sigma: \mathbb{N}^2 \to B(\clh)$ be an isometric representation, and let $V: [0,\infty)^2 \to B(L^2([0,1)^2) \bigotimes \clh)$ be the isometric representation induced by $\sigma$. Let $s,t \geq 0$. There exist $m,n \in \bbn$ such that $m \leq s < m+1$ and $ n \leq t < n+1$. Let 
	\begin{align*}
	R_{11} &:= [0,s-m) \times [0,t-n)\\ 
	R_{12}&: = [0,s-m) \times [t-n,1)\\
	 R_{21}&: = [s-m,1)\times [0,t-n)\\ 
	 R_{22}&=[s-m,1) \times [t-n,1).\end{align*} A short calculation, as in the case $d=1$, shows that for any $\xi \in L^2([0,1)^2) \bigotimes \clh$,

\[ V_{(s,t)}^* \xi (x,y) = \begin{cases}
	\sigma(m+1,n+1)^* \xi(x-s+m+1, y-t+n+1), ~\textrm{if~} (x,y) \in R_{11}, \\
	\sigma(m+1,n)^* \xi(x-s+m+1, y-t+n), ~ \textrm{if~}(x,y) \in R_{12} , \\
	\sigma(m,n+1)^* \xi(x-s+m, y-t+n+1), ~ \textrm{if~}(x,y) \in R_{21}, \\
	\sigma(m,n)^{*}\xi(x-s+m, y-t+n), ~ \textrm{if~}(x,y) \in R_{22}. \\
\end{cases} \]
		Define 
	\[\sigma^{(1)}: \mathbb{N} \to B(\clh), ~ \sigma^{(1)}(m):= \sigma(m,0), ~ m \in \mathbb{N}\]
	\[\sigma^{(2)}: \mathbb{N} \to B(\clh), ~ \sigma^{(2)}(n):= \sigma(0,n), ~ n \in \mathbb{N}\]
	Let $V^{(1)}$ and $V^{(2)}$ be the induced representations on $L^2([0,1)) \bigotimes \clh$ due to $\sigma^{(1)}$ and $\sigma^{(2)}$ respectively. Let $U: L^2([0,1)^2, \clh) \to L^2([0,1)^2, \clh)$ be  the flip defined by \[U \xi(x,y) = \xi(y,x)\] for all $\xi \in L^2([0,1)^2, \clh)$. Then,
	\begin{align*}
		U(1_{L^2([0,1))} \otimes V^{(1)}_s)U^* &= V_{(s,0)}, ~ \forall s \geq 0, \\
		1_{L^2([0,1))} \otimes V^{(2)}_t &= V_{(0,t)}, ~ \forall t \geq 0.
	\end{align*}
	
	\begin{ppsn} \label{2-d commutant}
		Suppose $\sigma$ is strongly pure, i.e, $\sigma^{(1)}$ and $\sigma^{(2)}$ are both pure isometric representations. Then,
		\[ \clm(V)' = 1_{L^2([0,1)^2)} \otimes \clm(\sigma)'\]
	\end{ppsn}
	
	\textit{Proof.} Since we have assumed  $\sigma^{(1)}$ and $\sigma^{(2)}$ are pure, it follows from   Prop. \ref{1-d commutant} that  $\clm(V^{(1)})' = 1_{L^2([0,1))} \otimes \clm(\sigma^{(1)})'$ and  $\clm(V^{(2)})' = 1_{L^2([0,1))} \otimes \clm(\sigma^{(2)})'$. Thus,
	\begin{align*}
		C^* \{ V_{(s,0)}| s \geq 0\}' &= U(B(L^2([0,1))) \otimes 1_{L^2([0,1))} \otimes \clm(\sigma^{(1)})' )U^* \\
		&= 1_{L^2([0,1))} \otimes B(L^2([0,1))) \otimes \clm(\sigma^{(1)})' \\
		C^* \{ V_{(0,t)}| t \geq 0\}' &= B(L^2([0,1))) \otimes 1_{L^2([0,1))} \otimes \clm(\sigma^{(2)})' .
	\end{align*}
	However, $\clm(V)' = C^* \{V_{(s,0)}| s \geq 0\}' \cap C^* \{V_{(0,t)}| t \geq 0\}'$. Thus, we have,
	\begin{align*}
		\clm(V)' &= 1_{L^2([0,1))} \otimes 1_{L^2([0,1))} \otimes (\clm(\sigma^{(1)})' \cap \clm(\sigma^{(2)})') \\ 
		&= 1_{L^2([0,1)^2)} \otimes \clm(\sigma)'.
	\end{align*} \hfill $\Box$
	
	\begin{lmma} \label{2-dim iff}
		Let $\eta: \mathbb{N}^2 \to \clh$ be a map. For $m \leq s< m+1$ and $n \leq t <n+1$, we define
		\[\xi_{(s,t)}^{\eta}(x,y): = \begin{cases} \eta_{(m,n)}, \textrm{ for } (x,y) \in [0,1-s+m) \times [0,1-t+n), \\
			\eta_{(m+1,n)}, \textrm{ for } (x,y) \in [1-s+m,1) \times [0,1-t+n), \\
			\eta_{(m,n+1)}, \textrm{ for } (x,y) \in [0,1-s+m)\times [1-t+n,1), \\
			\eta_{(m+1,n+1)}, \textrm{ for } (x,y) \in [1-s+m,1) \times [1-t+n,1). \end{cases} \]
		Then, $\xi^{\eta}= \{ \xi_{(s,t)}^{\eta}\}$ is an additive cocycle of $V$ iff $\eta$ is an additive cocycle of $\sigma$.
	\end{lmma}
	\textit{Proof.} Similar to Lemma \ref{iff}. \hfill $\Box$
	
	\begin{rmrk}
		Note that in order to define an element $\eta$ in $\cla(\sigma)$, it is sufficient to specify $\eta_{(1,0)} \in ker(\sigma(1,0)^*)$ and $\eta_{(0,1)} \in ker(\sigma(0,1)^*)$ such that \[\eta_{(1,0)} + \sigma(1,0) \eta_{(0,1)} = \eta_{(0,1)} + \sigma(0,1) \eta_{(1,0)}.\] 
	\end{rmrk}
	
	\begin{ppsn} \label{2-d additive cocycles}
		Let $\sigma^{(1)}$ and $\sigma^{(2)}$ be pure. The map
		\[\cla(\sigma) \ni \eta \mapsto \xi^{\eta} \in \cla(V)\]
		is an isomorphism.
	\end{ppsn}
	\textit{Proof.} The map is clearly injective. To see it is a bijection, let $\xi= \{ \xi_{(s,t)}\}_{s,t \geq 0}$ be an additive cocycle of $V$. Then, $\{ \xi_{(0,t)}\}_{t \geq 0}$ is an additive cocycle of  the $1$-parameter isometric representation $\{ V_{(0,t)}| t \geq 0\} = \{ 1_{L^2([0,1))} \otimes V^{(2)}_t | t \geq 0\}$. Thus, there exists $f \in L^2([0,1))$ and $\xi^{(2)} \in \cla(V^{(1)})$ such that $\xi_{(0,t)}=f\otimes \xi^{(2)}_t$, i.e. $\xi_{(0,t)}(x,y)=f(x) \xi^{(2)}_t(y)$. However, by Prop. \ref{1-d bijection}, there exists a unique $\eta^{(2)} \in \cla(\sigma^{(2)})$ such that $\xi^{(2)}= \xi^{(\eta^{(2)})}$. Therefore, for $n \leq t < n+1$,
	\[ \xi_{(0,t)}(x,y) = \begin{cases}
		f(x) \eta^{(2)}_n, ~\textrm{if ~} (x,y) \in [0,1) \times [0, 1-t+n), \\
		f(x) \eta^{(2)}_{n+1}, ~ \textrm{if ~}(x,y) \in [0,1) \times [1-t+n,1).
	\end{cases}\] 
	Similarly, $\{ \xi_{(s,0)}\}_{s \geq 0}$ is an additive cocycle of $\{ V_{(s,0)}\}_{s \geq 0} = \{ U(1_{L^2([0,1))} \otimes V^{(1)}_s)U^*\}_{ s \geq 0}$. This implies that $\{ U^* \xi_{(s,0)}\}_{s \geq 0}$ is an additive cocycle of $\{ 1_{L^2([0,1))} \otimes V^{(1)}_s | s \geq 0\}$. Thus, there exists $g \in L^2([0,1))$ and $\xi^{(1)} \in \cla(V^{(1)})$ such that $
		U^* \xi_{(s,0)}(x,y) = (g \otimes \xi^{(1)}_s) (x,y)$, and hence, $
		\xi_{(s,0)}(x,y) = g(y) \xi^{(1)}_s(x).$
	
	Again, by Prop. \ref{1-d bijection}, there exists $\eta^{(1)} \in \cla(\sigma^{(1)})$ such that, for $m \leq s < m+1$,
	\[ \xi_{(s,0)}(x,y) = \begin{cases}
		g(y) \eta^{(1)}_m, ~ \textrm{if~}(x,y) \in [0,1-s+m) \times [0,1), \\
		g(y) \eta^{(1)}_{m+1}, ~ \textrm{if ~}(x,y) \in [1-s+m) \times [0,1).
	\end{cases}\]
	For $0 \leq s,t <1$,
	\begin{align*}
		\xi_{(s,t)}(x,y) &= \xi_{(s,0)}(x,y) + V_{(s,0)}\xi_{(0,t)}(x,y)  \\
		&=
		\begin{cases}
			0, ~~~~~(x,y) \in [0,1-s) \times [0,1-t), \\
			g(y)\eta^{(1)}_1, ~ (x,y) \in [1-s,1) \times [0,1-t), \\
			f(x+s)\eta^{(2)}_1, ~ (x,y) \in [0,1-s) \times [1-t,1), \\
			g(y)\eta^{(1)}_1 + \sigma(1,0) f(x+s-1)\eta^{(2)}_1, ~ (x,y) \in [1-s,1) \times [1-t,1).
		\end{cases} \\
		\end{align*}
		Also, \begin{align*}
		 \xi_{(s,t)}(x,y) &= \xi_{(0,t)}(x,y) + V_{(0,t)}\xi_{(s,0)}(x,y) \\
		&=
		\begin{cases}
			0, ~ (x,y) \in [0,1-s) \times [0,1-t), \\
			g(y+t)\eta^{(1)}_1, ~ (x,y) \in [1-s,1) \times [0,1-t), \\
			f(x)\eta^{(2)}_1, ~ (x,y) \in [0,1-s) \times [1-t,1), \\
			f(x)\eta^{(2)}_1 + \sigma(0,1)g(y+t-1) \eta^{(1)}_1, ~ (x,y) \in [1-s,1) \times [1-t,1).
		\end{cases} 
	\end{align*}
	Combining both, we get, given $0 \leq t <1$,  $g(y)\eta^{(1)}_1 = g(y+t)\eta^{(1)}_1$, for all almost all $ y \in [0,1-t)$. Similarly, given $s \in [0,1)$, $f(x)\eta^{(2)}_1= f(x+s)\eta^{(2)}_1$ for almost all $x \in [0,1-s)$. Thus, $f\eta^{(2)}_1$ and $g\eta^{(1)}_1$ are both constant on $[0,1)$, i.e, there exist $c_0,d_0 \in \mathbb{C}$ such that $f(x)\eta^{(2)}_1 = d_0\eta^{(2)}_1$ for almost all $x \in [0,1)$, and $ g(y)\eta^{(1)}_1=c_0\eta^{(1)}_1$ for almost all $y \in [0,1)$. Hence, for $m \leq s <m+1$ and $n 
	\leq t < n+1$, we have
	\begin{align*}
		\xi_{(s,0)}(x,y) &= \begin{cases}
			c_0 \eta^{(1)}_m, ~ (x,y) \in [0,1-s+m) \times [0,1), \\
			c_0 \eta^{(1)}_{m+1}, ~ (x,y) \in [1-s+m,1) \times [0,1).
		\end{cases} \\
		\xi_{(0,t)}(x,y) &= \begin{cases}
			d_0 \eta^{(2)}_n, ~ (x,y) \in [0,1) \times [0,1-t+n), \\
			d_0 \eta^{(2)}_{n+1}, ~ (x,y) \in [0,1) \times [1-t+n,1).
		\end{cases} 
	\end{align*}	
	Define $\eta_{(m,n)} := c_0 \eta^{(1)}_m + \sigma(m,0) d_0 \eta^{(2)}_n$ for all $m,n \in \mathbb{N}$. Since $\eta^{(1)}_m \in ker(\sigma^{(1)*})$ and $\eta^{(2)}_n \in ker(\sigma^{(2)*})$, we have $\sigma(m,n)^* \eta_{(m,n)} =0$. Also,
	\begin{align*}
		\eta_{(1,0)} + \sigma(1,0) \eta_{(0,1)} &= c_0 \eta^{(1)}_1 + \sigma(1,0) d_0 \eta^{(2)}_0 + \sigma(1,0)( c_0 \eta^{(1)}_0 + \sigma(0,0) d_0 \eta^{(2)}_1) \\
		&= c_0 \eta^{(1)}_1 + \sigma(1,0)  d_0 \eta^{(2)}_1 \\
		&= \eta_{(1,1)}. \end{align*}
		Moreover,
				\begin{align*}
		\eta_{(0,1)} + \sigma(0,1) \eta_{(1,0)} &= c_0 \eta^{(1)}_0 + \sigma(0,0) d_0 \eta^{(2)}_1 + \sigma(0,1)( c_0 \eta^{(1)}_1 + \sigma(1,0) d_0 \eta^{(2)}_0) \\ 
		&= d_0 \eta^{(2)}_1 + \sigma(0,1) c_0 \eta^{(1)}_1 \\
		&= \xi_{(0,1)}(x,y) + \sigma(0,1) \xi_{(1,0)}(x,y) \\
		&= \xi_{(1,1)}(x,y) \\
		&=\xi_{(1,0)}(x,y)+\sigma(1,0)\xi_{(0,1)}(x,y)\\
		&=c_o\eta^{(1)}_1+\sigma(1,0)d_{0}\eta^{(2)}_{1}\\
		 &= \eta_{(1,1)}.
	\end{align*}
	Therefore, $\eta$ is an additive cocycle of $\sigma$. By Lemma \ref{2-dim iff}, $\xi^{\eta}$ is an additive cocycle of $V$. Since $\xi$ and $\xi^{\eta}$ are additive cocycles,  $\xi_{(s,0)}=\xi^{\eta}_{(s,0)}$ and $\xi_{(0,t)}=\xi^{\eta}_{(0,t)}$ for $s,t \geq 0$, it follows that $\xi=\xi^{\eta}$. Thus, the map is a bijection. That concludes the proof. \hfill $\Box$
	
	\begin{rmrk} \label{strongly pure}
		Recall that the representation $\sigma$ is called strongly pure if both $\sigma^{(1)}$ and $\sigma^{(2)}$ are pure isometric representations.. Suppose $\sigma$ is a pure isometric representation of $\bbn^2$ on $\clh$. Let $a,b \in \bbn^2$ be order units for $\bbz^2$, i.e., for every $x \in \mathbb{Z}^2$, there exist $m,n \in \mathbb{N}$ such that $ma - x$ and $nb - x$ both belong to $\mathbb{Z}^{2}$.  Define $\tilde{\sigma} : \bbn^2 \to B(\clh)$ by $\tilde{\sigma}(1,0) := \sigma(a)$ and $\tilde{\sigma}(0,1) := \sigma(b)$. Since $\sigma$ is a pure isometric representation and $a$ is an order unit, 
		\[\bigcap_{m \in \mathbb{N}} Ran(\tilde{\sigma}(m,0)) =\bigcap_{m \in \bbn}Ran(\sigma(ma))=\bigcap_{(m,n) \in \bbn^2}Ran(\sigma(m,n))=\{0\}.\]
		Therefore, $\tilde{\sigma}^{(1)}$ is a pure isometric representation of $\bbn$. Similarly, $\tilde{\sigma}^{(2)}$ is a pure isometric representation of $\bbn$. Thus, $\tilde{\sigma}$ is strongly pure. 
		
		 Choose order units such that  semigroup of $\bbz^2$ generated by $a$ and $b$ is contained in the semigroup $\bbn^2$ and also spans $\bbz^2$. For example, we can let $a=(1,1)$ and $b=(2,1)$. By Prop. \ref{smaller cone}, we have $dim(\cla(\sigma)) = dim(\cla(\tilde{\sigma}))$ and $\clm(\sigma)' = \clm(\tilde{\sigma})'$.
	\end{rmrk}

	\section{Examples}

	In this section, we prove Thm. \ref{main_intro}. First, we explain that Thm. \ref{main_intro} follows if we prove  the analogous two parameter discrete version. The reduction is explained below. 
	Let $P$ be a closed, convex cone in $\bbr^d$ that is spanning and pointed. Assume that $d \geq 2$. 
	\begin{enumerate}
		\item[(1)] Since $P$ is spanning and pointed, without loss of generality, we can assume that $P \subset \bbr_{+}^{d}$. Thanks to Prop. \ref{smaller cone}, it suffices to prove Thm. \ref{main_intro}
		when $P=\bbr_{+}^{d}$. Hereafter, assume that $P=\bbr_{+}^{d}$.
		
		\item[(2)] Suppose $V:\bbr_{+}^{2} \to B(\clh)$ is an isometric representation. Let $W:\bbr_{+}^{d} \to B(\clh)$ be defined by 
		\[
		W_{(t_1,t_2,\cdots,t_d)}:=V_{(t_1,t_2)}.\]
		Then, it is not difficult to show that $\mathcal{A}(W) \cong \mathcal{A}(V)
		$. Hence,  $Ind(W)=Ind(V)$. Clearly, $W$ is irreducible if and only if $V$ is irreducible. To denote the dependence of $W$ on $V$, we denote $W$ by $W^{V}$. 
		
		Moreover, for isometric representations $V_1$ and $V_2$  of $\bbr_{+}^{2}$,  $W^{V_1}$ and $W^{V_2}$ are unitarily equivalent if and only if $V_1$ and $V_2$ are unitarily equivalent.
		Thus, it suffices to prove Thm. \ref{main_intro} under the assumption that $P=\bbr_{+}^{2}$.
		
		\item[(3)] Thanks to Prop. \ref{2-d additive cocycles} and Prop. \ref{2-d commutant}, to prove Thm. \ref{main_intro} when $P=\bbr_{+}^{2}$, it suffices to produce, for any given $k$, a continuum of strongly pure irreducible 
		isometric representations of $\bbn^2$ with index $k$. Remark \ref{strongly pure} allows us to drop the requirement that the desired irreducible isometric representations of $\bbn^2$ need to be strongly pure. 
	\end{enumerate}

	With the discussion above, the problem now boils down to finding pure isometric representations of $\bbn^2$ that are irreducible and whose space of additive cocycles have dimension $k$ for $k \in \{0,1,2,3,\cdots,\} \cup \{\infty\}$.

	\begin{ppsn}\label{example}
		Let $\clh$ be a separable Hilbert space. Suppose $d \in \{1,2,\cdots \} \cup \{ \infty\}$. Let $\{P_i| 1 \leq i \leq d\}$ be a family of mutually orthogonal projections on $\clh$ such that $\displaystyle \sum_{i=1}^{d}P_i =1$. For each $k \geq 1$, define $\displaystyle Q_k := 1 - \sum_{i=1}^k P_i$. Let $U$ be a unitary on $\clh$. Define an isometric representation $\sigma: \bbn^2 \to B(\clh \bigotimes \ell^2(\bbn))$ by
		\[ \sigma(1,0):= 1 \otimes S, ~~ \sigma(0,1):= \sum_{i=1}^{d} UP_i \otimes S^{i-1}, \]
		where $S$ is the usual shift operator on $\ell^2(\bbn)$. Then,
		\begin{enumerate}
			\item $dim(\cla(\sigma)) = dim(\{ x \in \clh| x \in ker(U-1) \textrm{ and } \sum_{i=1}^{d}||Q_i x||^2 < \infty \})$, and 
			\item $\clm(\sigma)' = C^*(\{U,P_i| 1 \leq i \leq d\})' \otimes 1$.
		\end{enumerate}
	\end{ppsn}	
	
	\textit{Proof.} Let $\{ \delta_i\}_{i \geq 0}$ be the standard orthonormal basis for $\ell^2(\bbn)$. Let $\xi = \{ \xi_{(m,n)}\}_{(m,n) \in \bbn^2}$ be an additive cocycle of $\sigma$. Since $\xi_{(1,0)} \in ker(\sigma(1,0)^{*})$, $\xi_{(1,0)} = x \otimes \delta_0$, for some $x \in \clh$. Let $\xi_{(0,1)} = \displaystyle \sum_{j \geq 0} y_j \otimes \delta_j$. Since $\sigma(0,1)^{*}\xi_{(0,1)}=0$,  
	\[ \sum_{k \geq 0}( \sum_{j \geq k} P_{j-k+1} U^* y_j) \otimes \delta_k =0.\]
	This gives us that for each $j \geq 0$, 
	\[ U^* y_j \in  ker( \sum_{i=1}^{ j+1}P_i) = Ran(Q_{j+1}).\]
	Now, since $\xi$ is an additive cocycle, it satisfies
	\[ \xi_{(1,0)} + \sigma(1,0) \xi_{(0,1)} = \xi_{(0,1)} + \sigma(0,1) \xi_{(1,0)},\]
	which in turn implies 
	\[ U^* y_j = U^* x - (\sum_{i=1}^{j+1}P_{i}) x, \forall j \geq 0.\]
	Since $U^{*}y_j \in Ran(Q_{j+1})$,    we get $P_i x = P_i U^* x$ for every $i$,  and consequently, for every $j$, $U^* y_j = Q_{j+1}U^*x$, i.e $y_j=UQ_{j+1}U^{*}x$. Since $P_i(x-U^*x) =0 $, for all $i$, and $\sum_{i}P_i=1$, we have $U^*x=x$. Hence, $y_j=UQ_{j+1}x$. The fact that $\sum_{j \geq 0} y_j \otimes \delta_j \in \clh \bigotimes \ell^2(\bbn)$ implies $\sum_{i=1}^{d}||Q_i x||^2 < \infty$.
	
	Conversely, choose $x \in ker(U-1)$ such that $\sum_{i=1}^{d}||Q_i x||^2 < \infty$. Let $\eta_{(1,0)}= x \otimes \delta_0$ and $\eta_{(0,1)}= \sum_{j \geq 0} UQ_{j+1}x \otimes \delta_j$. It is routine to check that 
		\[ \eta_{(1,0)} + \sigma(1,0) \eta_{(0,1)} = \eta_{(0,1)} + \sigma(0,1) \eta_{(1,0)}.\]
	Thus, there exists an additive cocycle $\xi:=\{\xi_{(m,n)}\}_{(m,n) \in \bbn^2}$ such that $\xi_{(1,0)}=x \otimes \delta_0$ and $\xi_{(0,1)}=\sum_{j \geq 0}UQ_{j+1}x \otimes \delta_j$. 
	
	Define a map $\{ x \in \clh| x \in ker(U-1) \textrm{ and } \sum_{i=1}^{d}||Q_i x||^2 < \infty \} \to \cla(\sigma)$ by
	\[ x \mapsto \xi^x,\]
	where $\xi^x_{(1,0)}:= x \otimes \delta_0$ and $\xi^x_{(0,1)}:=  \displaystyle \sum_{j \geq 0} UQ_{j+1}x \otimes \delta_j$. It is now obvious that this map is an isomorphism. Therefore, \[dim(\cla(\sigma)) = dim(\{ x \in \clh| x \in ker(U-1) \textrm{ and } \sum_{i=1}^{d}||Q_i x||^2 < \infty \}).\] That concludes the first part of the proof.
	
	Let $T \in \clm(\sigma)'$. Since, $T \in C^*\{ \sigma(1,0)=1 \otimes S\}'$, $T$ is of the form $T=T_0 \otimes 1$, for some $T_0 \in B(\clh)$. Also, $T\sigma(0,1) = \sigma(0,1)T$. Thus,
	\[ \sum_{i \geq 1}T_0UP_i \otimes S^{i-1} = \sum_{i\geq 1}UP_iT_0 \otimes S^{i-1}\]
	Hence, $T_0UP_i= UP_iT_0$, for all $i \geq 1$. Summing over all $i$, we get
	\[T_0 U= UT_0,\]
	and that implies, for all $i \geq 1$,
	\[T_0 P_i = P_iT_0.\]
	Thus, $\clm (\sigma)' \subseteq C^*(\{U,P_i| 1 \leq i \leq d\})' \otimes 1$. Clearly, the reverse side of the inclusion holds as well, and we get 
	\[ \clm(\sigma)' = C^*(\{U,P_i| 1 \leq i \leq d\})' \otimes 1.\] \hfill $\Box$	
	
	\begin{rmrk}
		Note that, in Prop. \ref{example}, when $\clh$ is finite-dimensional, \[
		Ind(\sigma)=dim(\cla(\sigma)) = dim(ker(U-1)).\] 
	\end{rmrk}

	We use  Prop. \ref{example} to produce concrete examples of irreducible isometric representations of $\bbn^2$ with any given index. 
	
	\textbf{Example 1:} 	For this example,  we refer to the work of Albeverio and Rabanovich \cite{ALBEVERIO201935}. Let $B_3$ denote Artin's braid group,
	\[B_3= \langle S,J| S^2= J^3 \rangle\]
	
	Let $m$ be a positive integer. Theorem 5 of \cite{ALBEVERIO201935} asserts that  there exists a non-empty open set $\Omega$ in a Euclidean space and a family  of irreducible unitary representations $\{\pi_h\}_{h \in \Omega}$ of $B_3$ on $\bbc^{6m}$ such that 
	\begin{enumerate}
	\item[(1)] for $h \neq k$, $\pi_h$ and $\pi_k$ are not unitarily equivalent, and
	\item[(2)] for $h \in \Omega$, $dim (ker(\pi_h(J)-1))=2m$.
	\end{enumerate}
	For the explicit expression of the representation $\pi_h$, the reader is referred to Section 3 of  \cite{ALBEVERIO201935}. For $h \in \Omega$, set $P_h:=\frac{1+\pi_h(S)}{2}$, and $U_h:=\pi_h(J)$, and define an isometric representation $\sigma_h:\bbn^2 \to B(\bbc^{6m} \otimes \ell^2(\bbn))$ by setting
	\[
	\sigma_h(1,0)=1\otimes S; ~~\sigma_h(0,1)=U_hP_h \otimes 1+U_h(1-P_h)\otimes S.\]
		Let $h \in \Omega$. It follows from $(2)$ and from Prop. \ref{example} that $Ind(\sigma_h)=2m$ for every $h \in \Omega$. By Prop. \ref{example} and the fact that $\pi_h$ is irreducible, it follows that $\sigma_h$ is irreducible. 
	
	Let $h,k \in \Omega$ be given. Suppose $T$ is a unitary operator that intertwines $\sigma_h$ and $\sigma_k$, i.e  $T\sigma_h(\cdot)=\sigma_k(\cdot)T$. Then, $T(1\otimes S)=(1\otimes S)T$. Hence, $T$ is of the form $T=T_0 \otimes 1$ for some unitary operator $T_0$ on $\bbc^{2m}$. The equality 
	\[
	(T_0 \otimes 1)\sigma_h(0,1)=\sigma_k(0,1)(T_0 \otimes 1)\] leads to the conclusion $T_0U_h=U_kT_0$ and $T_0P_h=P_kT_0$. In other words, $\pi_h$ and $\pi_k$ are unitarily equivalent. Thus, the isometric representations $\sigma_h$ and $\sigma_k$ are not unitarily equivalent if $h \neq k$. 	 This gives a plethora of irreducible examples with even index that are not unitarily equivalent.

	It is not difficult to construct other examples as the following two classes of examples show. 
	
		\textbf{Example 2: }
			 Let $\clh$ be a finite-dimensional Hilbert space with an orthonormal basis $\{e_1, e_2, ..., e_n\}$. Let $P_i \in B(\clh)$ be the projection onto the subspace spanned by $e_i$, for $i=1,2,...,n$. Choose $a \in \clh$ such that $\langle a| e_i \rangle \neq 0$ for every $i = 1,2,...,n$. Let $P_a \in B(\clh)$ be the projection onto the subspace spanned by $a$,  i.e. $P_a(x) := \langle x | a \rangle a $, for all $x \in \clh$. Define $U_a \in B(\clh)$ by $U_a= 1 - 2 P_a$. Then, $U_a$ is a self-adjoint unitary on $\clh$ and $dim(ker(U_a-1)) = dim(ker(P_a)) = n-1$. Also, $C^*(\{U_a, P_1, ..., P_n\})' = \bbc$. Define an isometric representation $\sigma^{(a)}:\bbn^2 \to B(\clh \otimes \ell^2(\bbn))$ by  
			 \[ \sigma^{(a)}(1,0) := 1 \otimes S, ~~ \sigma^{(a)}(0,1) := \sum_{i=1}^{n} U_a P_i \otimes S^{i-1}. \]			 
	It follows from Prop. \ref{example} that $Ind(\sigma^{(a)})=n-1$, and $\sigma^{(a)}$ is irreducible 

Let $a, b \in \clh$ be such that $\langle a|e_i \rangle \neq 0$, $\langle b|e_i \rangle \neq 0$, for all $1 \leq i \leq n$, and there exists $k \in \bbn$, $1 \leq k \leq n$ such that $| \langle a|e_k	\rangle | \neq | \langle b|e_k \rangle |$. 
We claim that $\sigma^{(a)}$ and $\sigma^{(b)}$ are not unitarily equivalent. Let us assume there exists a unitary $T: \clh \bigotimes \ell^2(\bbn) \to \clh \bigotimes \ell^2(\bbn)$ such that 
\[ T\sigma^{(a)}(m,n) = \sigma^{(b)}(m,n)T\]
for $(m,n) \in \bbn^2$. 

Since $T$ commutes with $\sigma^{(a)}(1,0) = \sigma^{(b)}(1,0) = 1 \otimes S$, it follows that $T= T_0 \otimes 1$, for some unitary $T_0 \in B(\clh)$. Also, $(T_0 \otimes 1)\sigma^{(a)}(0,1) = \sigma^{(b)}(0,1)(T_0 \otimes 1)$ gives
\begin{equation}
\label{1234}T_0 U_a P_i = U_b P_i T_0, \textrm{ for all } 1 \leq i \leq n.\end{equation}
Adding them, we get 
\begin{equation}
\label{5678}T_0 U_a = U_b T_0.\end{equation}
Eq. \ref{5678} and Eq. \ref{1234} imply that $T_0P_i=P_iT_0$ for every $i$. Thus, \[T_0 e_i = \lambda_i e_i\]
for some $\lambda_i  \in \bbt$, for all $1 \leq i \leq n$. Eq. \ref{5678} implies $T_0 P_a e_k = P_b T_0 e_k$ for all $k$, which simplifies to 
\[ \lambda_i \langle e_k| a \rangle \langle a| e_i \rangle = \lambda_k \langle e_k|b \rangle \langle b|e_i \rangle, ~ \forall i\]
Thus, we get  $| \langle a|e_k \rangle | = | \langle b|e_k \rangle|$ for all $k$,  which is a contradiction. Therefore, $\sigma^{(a)}$ and $\sigma^{(b)}$ are not unitarily equivalent whenever  there exists $k$ such that $| \langle a|e_k \rangle| \neq | \langle b|e_k \rangle|$. 

\textbf{Example 3:} Let $\clh$ be an infinite dimensional Hilbert space with an orthonormal basis $\{e_n\}_{n \in \bbn}$. Define $P_i: \clh \to \clh$ by $P_i e_j = \delta_{ij} e_i$, for all $i \in \bbn$. Choose any $a \in \clh$ with $\langle a| e_i \rangle \neq 0$ for all $i \in \bbn$, and let $P_a e_i := \langle e_i|a \rangle a$, for all $i$. Set $U_a= 1-2P_a$. Define an isometric representation $\sigma^{(a)}: \bbn^2 \to B(\clh \bigotimes \ell^2(\bbn))$ by
\[ \sigma^{(a)}(1,0):= 1 \otimes S, ~~ \sigma^{(a)}(0,1):= \sum_{i=1}^{\infty} U_aP_i \otimes S^{i-1}.\]
Note that $C^{*}\{U_a, P_i:i=1,2,\cdots\}^{'}=\bbc$. Hence, by Prop. \ref{example}  $\sigma^{(a)}$ is an irreducible isometric representation of $\bbn^2$. Again by Prop. \ref{example}, $\sigma^{(a)}$ has index equal to the dimension of the space $\displaystyle \Big \{ x \in ker(U_a -1)| \sum_{n\geq 2} (n-1)||P_n x||^2 < \infty \Big \}$ and the latter has infinite dimension. Thus, $Ind(\sigma^{(a)})=\infty$.  Also, just like in the finite-dimensional case, if we choose $a,b \in \clh$ such that for some $k \in \bbn$, $| \langle a| e_k \rangle | \neq | \langle b| e_k \rangle |$, then $\sigma^{(a)}$ and $\sigma^{(b)}$ will not be unitarily equivalent.

We encompass all of the above in the following theorem.
	
\begin{thm}
\label{last theorem}
For each $k \in \{0, 1, 2,\cdots,\} \cup \{\infty\}$, there is a continuum of irreducible isometric representations of $\bbn^2$ that has index $k$. 
\end{thm}	

Note that Thm. \ref{main_intro} is now immediate from Thm. \ref{last theorem} and the discussions made at the beginning of this section.  Thm. \ref{main_intro} together with Remark \ref{on ccr} imply Thm. \ref{main prime}.

	\bibliography{references_piyasa}
	\bibliographystyle{amsplain}

\end{document}